# Adaptive Bayesian Optimization Algorithm for Unpredictable Business Environments.

Sarit Maitra
Sarit.maitra@gmail.com
Alliance Business School, Alliance University
Bangalore, India

**Abstract**

This paper presents an innovative optimization framework and algorithm based on the Bayes theorem, featuring adaptive conditioning and jitter. The adaptive conditioning function dynamically modifies the mean objective function in each iteration, enhancing its adaptability. The mean function, representing the model's best estimate of the optimal value for the true objective function, is adjusted based on observed data. The framework also incorporates an adaptive acquisition jitter function, enhancing adaptability by adjusting the jitter of the acquisition function. It also introduces a robust objective function with a penalty term, aiming to generate robust solutions under uncertainty. The evaluation of the framework includes single-objective, decoupled multi-objective, and combined multi-objective functions. Statistical analyses, including t-statistics, p-values, and effect size measures, highlight the superiority of the proposed framework over the original Bayes optimization. The adaptive nature of the conditioning function allows the algorithm to seamlessly incorporate new data, making it particularly beneficial in dynamic optimization scenarios.

## 1. INTRODUCTION

Global optimization is to find the best optimal values of an objective function over a feasible space of input variables. It identifies the global optimum, which is the set of input values that leads to the lowest value of the objective function. The local optimization finds the best solution near a starting point, while the global optimization aims to explore the entire feasible region. This is to ensure that the best solution is found across the search space. This makes global optimization problems often challenging because the objective function has multiple local optima, and traditional optimization methods often get stuck in a suboptimal solution. In the context of supply chain inventory policy development, this challenge is aggravated by the inherent complex structures along with the rising stochastic demand patterns, lead times, product variety, and shortened lifecycles. Traditional optimization methods, rooted in classical inventory management theory (Hameri & Paatela, 2005; Gupta & Maranas, 2003), often fall short of dealing with such challenges.

Inventory management is an important aspect of the supply chain, where simulation is often used to overcome the inflexible nature and find out the probabilistic measures (Thevenin et al., 2021; Wu & Frazier, 2019; Gruler et al., 2018; Azadi et al., 2019). Although historically, inventory management has primarily focused on deterministic demand models in a supply chain context (e.g., Babiloni & Guijarro, 2020; Stopková et al., 2019; Tamjidzad & Mirmohammadi, 2017 & 2015; Choi, 2013; etc.), in recent times, stochastic demand modeling has garnered serious attention from practitioners and academicians due to the increasing uncertainty in demand patterns. Thus, a shift can be observed from deterministic to probabilistic modeling and inference (e.g., Heikkinen et al., 2022; Li et al., 2021; Wang et al., 2020; Gholami et al., 2021; Pearce et al., 2022). The stochastic model inherently induces uncertainties and variations, which result in several outcomes with varying probabilities (Hasan et al., 2019; Sakki et al., 2022), as outlined in the simulation method. Though simulations are useful for simulating the complexities of real-world business scenarios, they are not very good at identifying the best decision variables (Kiuchi et al., 2020). Moreover, simulating real-world complexities requires substantial computational effort. Empirical studies show that the benefits of a simulation study can be realized when integrated with an effective optimization framework (Kastner et al., 2021; Huang et al., 2019).

The objective function of the optimization problem is often unique since it is an unknown analytic form, or a black box. This automatically reduces the number of optimization algorithms available for solving the problem. Moreover, to make this cost-effective, it is also preferable to get optimality with the fewest possible function evaluations (Roman et al., 2019). Bayes Optimization (BayesOpt) (Garnett et al., 2010) has emerged as a robust and effective optimization strategy, garnering attention from several researchers (e.g., Hosseini et al., 2020; Oyewola et al., 2022; Hosseini & Ivanov, 2020 & 2022, etc.). It provides probabilistic models with unique methodological qualities for modeling dependencies in complex networks (Hosseini & Ivanov, 2020). In the optimization space, numerous fields of study have experimented successfully with BayesOpt, such as the social sciences, ecology (Bac Dang et al., 2019; Lau et al., 2017); medicine (Haddawy et al., 2018; Wang et al., 2019); services and banking (De Sa et al., 2018; Tavana et al., 2018); the energy, defense, and robotics industries (Boutselis and McNaught, 2019; Munya et al., 2015); etc. Compared with other optimization frameworks, BayesOpt often significantly reduces the number of function evaluations (Kiuchi et al., 2020) and excels at modeling unpredictability and capturing non-linear causal links, enabling the drawing of conclusions from imprecise, uncertain, and incomplete data (Wan et al., 2019).

The fundamental principle of BayesOpt involves the construction of a surrogate function to represent the objective function. This is then explored using an acquisition function such as Expected Improvement (EI) or Lower Confidence Bound with Entropy Weighting (LCB) (Garnett, 2023; Hennig et al., 2022). Its great data efficiency enables efficient exploration and exploitation, eliminating the need for several evaluations (Wu et al., 2019). However, empirical studies underscored the limitations associated with the feasible initial point in BayesOpt strategies, such as Expected Improvement (EI) with constraints (Snoek, 2013; Gelbart et al., 2014). To address these challenges, researchers have already proposed changes such as integrated conditional expected improvement (Bernardo et al., 2011) and expected volume reduction (Picheny, 2014). These approaches aim to enhance the flexibility and applicability of BayesOpt by mitigating the constraints associated with the initial point requirement, thereby broadening the scope of problems that can be effectively addressed. However, finding feasible points consumes computational budget, potentially slowing the optimization process. BayesOpt also requires joint evaluation of objective function and constraints at candidate points to quantify its utility for finding a global optimum, which can be time-consuming (Ariafar et al., 2019).

Although BayesOpt has progressed from a new idea to a growing research field when it comes to supply chain resilience and risk assessment, there is not much research currently available when it comes to using inventory or supply chain optimization challenges. However, its limitations make its application and adaptation problematic in real-life scenarios. This work aims to introduce an adaptable BayesOpt framework that addresses practical issues and improves adaptation in dynamic business contexts. The adaptive framework addresses disconnected problems, non-closed form acquisition functions, and approximation concerns, addressing coordination issues and global optimality challenges. Non-closed form acquisition functions require approximations or numerical methods for optimal points, which can be computationally challenging and require advanced techniques. The accuracy of the surrogate model can be improved using sophisticated modeling approaches or model adaptation during the optimization phase.

This work introduces an adaptive framework that enhances algorithm flexibility and adaptability to changes in the optimization landscape. The framework dynamically modifies the mean function based on observed data, ensuring a balance between exploration and exploitation. The evaluation is conducted using statistical analyses and analytical measures, assessing the differences between the proposed framework and traditional approaches. The framework demonstrates efficiency in handling both decoupled and disconnected evaluations, a critical aspect in real-world optimization scenarios.

## 2. METHODOLOGY

BayesOpt is used in real-world scenarios like complex simulations, inventory management, engineering designs, and scientific research due to the uncertainty of the true objective function. Computed analytical formulations can be expensive or unavailable, making it difficult to measure goal function values at every input point due to budgetary, scheduling, or practical constraints.

## 2.1. Framework formulation

This framework employs BayesOpt, Gaussian Process Regression, and Expected Improvement for system optimization, with a dynamic conditioning mechanism that adjusts the mean function based on observed data, making it effective in uncertain or dynamically shifting optimization landscapes.

The objective function of the GPR model:

$$f(x) \sim GP(m(x), k(x,x')) \quad (1)$$

where $m(x)$ and $k(x,x')$ are the mean and covariance functions and denoted as:

$$\left.\begin{array}{l} m(x) = (f(x)) \\ k(x,x') = E\{[f(x) - m(x)][f(x') - m(x')]^T\} \end{array}\right] \quad (2)$$

$E[f(x)]$ is the expected value (average value) of the function $f(x)$ at the input $x$. Similarly, $E[(f(x) - m(x))(f(x') - m(x'))^T]$ is the expected value of the outer product of the deviations of $f(x)$ and $f(x')$ from their respective mean functions. When modeling the relationships between input data points, the choice of kernel, or covariance function, is essential. Matern52 is the suggested kernel since it allows for improved modeling of the underlying stochastic goal function.

$$k_{Matern52}(r) = \left(1 + \frac{\sqrt{5}r}{l} + \frac{\sqrt{5}r^2}{3l^2}\right) exp\left(-\frac{\sqrt{5}r}{l}\right) \quad (3)$$

where $r$ is the pairwise distance between input points and $l$ is a lengthscale parameter. This kernel has a smooth behavior and can capture complex patterns in the data. $\sqrt{5}$ is a scaling factor in the expression. It is a design choice that helps balance the trade-off between capturing complex patterns in the data and maintaining smoothness in the predictions. However, adaptive kernel selection (e.g., Roman et al., 2019) can be opted here to explore further.

The mean function captures the expected trend, and the covariance function models the relationships and uncertainties between different points. In the probabilistic context, the expectation operator provides a way to compute the average behavior of a random variable over all outcomes, weighted by their respective probabilities.

Robust objective function:

$$f_{robust}(x, iteration) = -mean(max(f(x) - 2, 0)) - 0.1 * mean(max(2 - f(x), 0)) \quad (4)$$

Here, the penalty term is introduced to the original objective function $f(x)$ to ensure robustness. The mean and max functions operate elementwise on the samples from the normal distribution, and the penalty term penalizes deviations below a threshold.

BayesOpt aims to find out the optimal input $X^*$ to minimize or maximize the objective function.

$$X^* = \arg\max_x f(x) \quad (5)$$

The introduction of an adaptive conditioning function that dynamically modifies the mean function based on observed data is a unique feature (Eq. 6). The optimization technique becomes much more flexible as additional data comes in (Eq. 7). It modifies the mean function $m(x)$ based on observed data:

$$m_c(x) = m(x) + \kappa(x) \quad (6)$$

here $m_c(x)$ = conditioned mean function. The conditioning function is:

$$\kappa(x) = \alpha * u(x) \quad (7)$$

where $\alpha$ = scaling factor and $u(x)$ captures the influence of observed data on the mean function. The adaptive approach dynamically adjusts parameters over iterations: $\kappa(x)$ (Eq. 8) and dynamic_jitter (Eq. 9).

$$\kappa(x) = \frac{constant\_value}{1 + 0.1 * iteration} \quad (8)$$

$$dynamic\_jitter = \frac{0.01}{1 + 0.1 * iteration} \quad (9)$$

To optimize the original objective function f(x), the combined function (Eq 10) seeks an input that balances prior expectations (represented by m(x)) with observed data (captured by κ(x)).

$$f_{combined}(x) = m_c(x) + f(x) \quad (10)$$

The objective is to attain the optimal input $x^*$ that maximizes or minimizes the combined objective function ($x^* = \arg max_x f_{combined}(x)$). This optimization aims to find an input that balances prior expectations (captured by $m(x)$) with observed data (captured by $\kappa(x)$), while optimizing the original objective function $f(x)$. This strategy enables the optimization algorithm to adjust its expectations as data accumulates, resulting in more informed conclusions over iterations.

EI function optimizes the combined objective function from Eq. 6:

$$EI(x) = \int_{-\infty}^{f_{min}} (f_{min} - f_{combined}(x))^+ \, p(f_{combined}(x)|D) \, df_{combined}(x) \quad (11)$$

Here, $f_{min}$ is the minimum observed value of the objective function, $p(f_{combined}(x)|D)$ is the posterior probability distribution of $f_{combined}(x)$ given the observed data D, and $(a)^+$ denotes the positive part of the function.

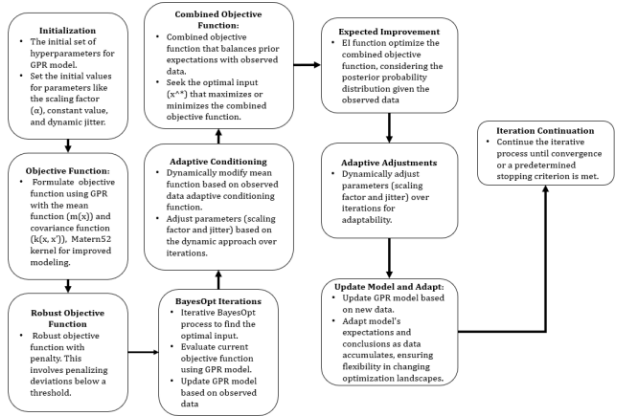

The optimization process is iterative, updating the model based on new data. The conditioning function ensures adaptability by dynamically adjusting parameters like scaling factor and jitter over iterations. This framework distinguishes itself from traditional BayesOpt, providing more informed conclusions over iterations. A simplified workflow that explains the procedures is shown in Fig. 1.

Figure 1: Single objective performance comparison

## 2.2. Analytical functions

Table 1 provides the pseudocode, which outlines the main steps and functions used. This function takes a set of hyperparameters (variance and lengthscale) and performs BayeOpt on the stochastic objective function with a GPR model using the Matern52 kernel. The goal is to find the hyperparameters that maximize profit by optimizing the stochastic objective function. All unknown quantities are represented as random variables in Bayesian inference. This is a powerful assumption because it permits views about these quantities to be represented by probability distributions reflecting their values. Inference then takes the form of an inductive process where the values are improved iteratively by considering observed facts and appealing to probabilistic identities..

This work introduces the robust objective function, which introduces a penalty term (penalty) based on the mean of the maximum of zero and the differences between a threshold (2) and the demand samples. The penalty term scales the factor (0.1) to emphasize the importance of robustness. The overall profit (assuming profit maximization is the objective) is the negative mean of the maximum of zero and the differences between the demand samples and the threshold, minus the penalty. This approach encourages the optimization algorithm to find solutions that not only optimize the original objective function but also penalize solutions that are resistant to uncertainties, promoting a robust solution over iterations.

ALGORITHM 1: Pseudocode - kernel hyperparameter tuning

1. **# Stochastic objective function**
2. def stochastic_objective_function(args):
3.    x = args[0][0]
4.    demand_samples = generate_demand_samples(x)
5.    profit = calculate_profit(demand_samples)
6.    return profit
7. **# Hyperparameter tuning objective**
8. def hyperparameter_tuning_objective(params):

9.  variance, lengthscale = params[0][0], params[0][1]
10. # BayesOpt
11. optimizer = BayesOpt(stochastic_objective_function, variance, lengthscale)
12. optimizer.run_optimization(iterations = 30)
13. return -get_best_objective_value(optimizer)
14. # Hyperparameter tuning
15. hyperparameter_optimizer = initialize_hyperparameter_optimizer(hyperparameter_tuning_objective)
16. hyperparameter_optimizer.run_optimization(max_iter=30)
17. best_hyperparameters = get_best_hyperparameters(hyperparameter_optimizer)
18. # variance and lengthscale
19. {best_hyperparameters[0]}, {best_hyperparameters[1]}')

The pseudocode (Algorithm 2) involves random samples from a normal distribution that is non-deterministic. This was intentional, keeping in mind that a deterministic approach may limit the exploration of the algorithm. Incorporating conditioning, adaptive expected improvement (EI), and robust BayesOpt aim at enhancing the optimization process in different ways.

ALGORITHM 2: Pseudocode - Robust objective function

1.  def robust_objective_function(args, iteration):
2.  x = args[0][0]
3.  demand_samples = np.random.normal(loc = x, scale = 1.0, size = 1000)
4.  # Robustness penalty term
5.  penalty = 0.1 * np.mean(np.maximum(2 - demand_samples, 0))
6.  # Original objective function
7.  profit = -np.mean(np.maximum(demand_samples - 2, 0))
8.  # Combined the above two functions
9.  combined_objective = profit - penalty
10. return combined_objective

Algorithm 3 displays the pseudocode for adaptive conditioning and the adaptive acquisition jitter function. This function dynamically adjusts the conditioning applied to the mean function of the GPR model based on the number of optimization iterations. It helps the algorithm adapt to exploration and exploitation as more data is accumulated, potentially leading to informed decisions over successive iterations. The jitter function dynamically adjusts the acquisition jitter used in the acquisition function. It controls the balance between exploration and exploitation.

ALGORITHM 3: Pseudocode - Adaptive conditioning

1.  function adaptive_conditioning(mu, constant_value, iteration):
2.  dynamic_constant = constant_value / (1 + 0.1 * iteration)
3.  modified_mean = mu + dynamic_constant
4.  return modified_mean
5.  function adaptive_acquisition_jitter(iteration):
6.  dynamic_jitter = 0.01 / (1 + 0.1 * iteration)
7.  return dynamic_jitter
8.  function main_optimization():
9.  # Initialization
10. num_iterations = 100
11. constant_value = 0.1
12. mu_initial = initial_mean()

13. # Optimization loop
14. for iteration in range(num_iterations):
15.     # Evaluate objective function and update observed data
16.     # GPR
17.     mu_current = update_mean(mu_initial)
18.     kernel = update_kernel()
19.     # Conditioning
20.     modified_mean = adaptive_conditioning(mu_current, constant_value, iteration)
21.     # BayesOpt using the modified mean and updated kernel
22.     # Acquisition function calculation and optimization step
23.     # Update observed data and GPR model
24. # End of the optimization loop

Parallelizing BayesOpt is not straightforward due to the inherent sequential nature of the optimization process, where each iteration depends on the results of the previous one. However, for Algorithm 3, parallelization can be achieved by running multiple independent optimizations concurrently, each with its instance of BayesOpt.

## 3. EMPIRICAL ANALYSIS & DISCUSSIONS

The study compares the BayesOpt process with added functions, focusing on the impact of adaptive conditioning and robustness considerations. Statistical tests assess the difference in results, while visualizations provide insights into convergence behavior and improvement rates. Figs. 2, 3, and 4 display the convergence and performance for single-objective, multi-objective, and decoupled multi-stochastic objectives consecutively. The system specifications used for this work: Google Cloud, Python v3.10.12, GPyOpt, Windows PC, Intel (R) Core (TM) i5-4570 CPU @ 3.20GHz, 16.0. GB.

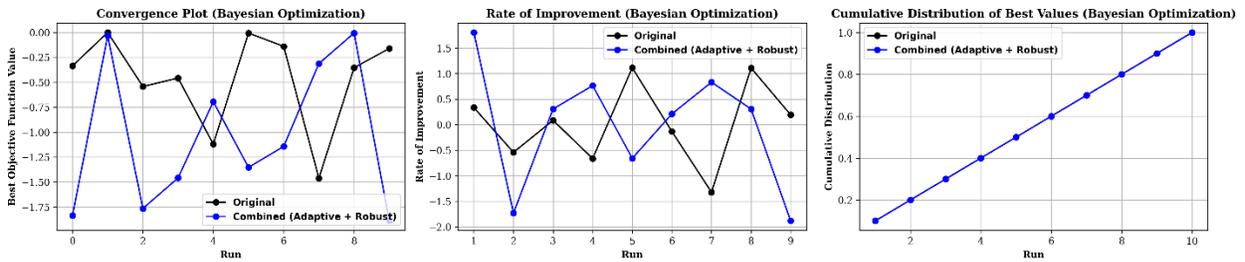

Figure 2: Single objective performance comparison

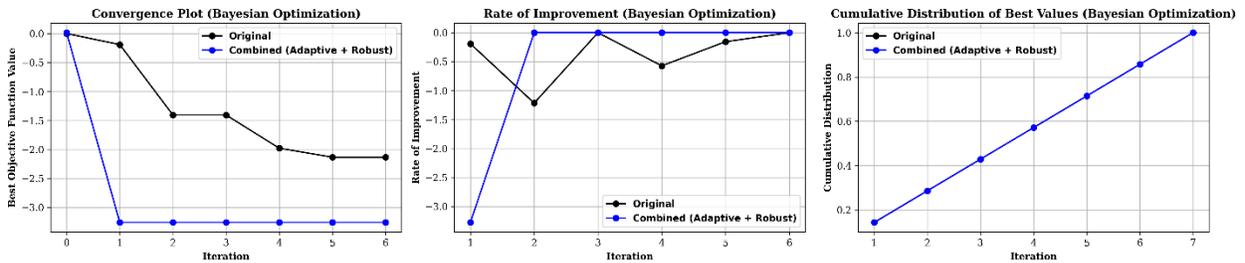

Figure 3: Multi-objective performance comparison

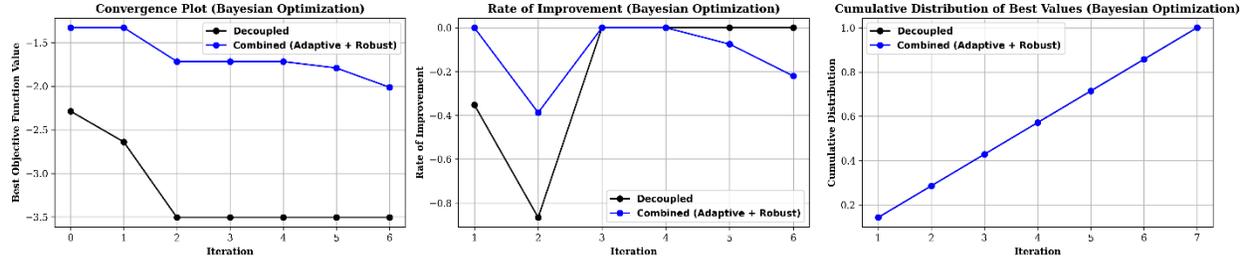

Figure 4: Multi-objective Decoupled BayesOpt performance comparison

The plots show no meaningful difference for a single objective function. However, compared to the original BayesOpt, we see a considerable improvement (Figs.3 and 4) in the convergence rate with the proposed changes for the multi-objective function. Higher values (middle plot) indicate faster convergence. The first plot (left side) displays a multi-objective case scenario where we observe a sharp drop indicating a significant improvement in the objective function value. This drop signifies that the algorithm has found a point in the search space that yields a much better (lower or higher, depending on the optimization goal) objective function value compared to the previous iterations. The analytical summary displayed in Table 4 compares the performances of each.

Table 1: Analytical summary of the optimization process.

| Optimization process | Best values | Rate of improvement | Cumulative distribution |
|---|---|---|---|
| Original BayesOpt- Single objective | -0.3345, 0.0, -0.5418, -0.4571, -1.1203, -0.0057, -0.1393, -1.4621, -0.3531, -0.1603 | -1.8375, -0.0322, -1.7657, -1.4591, -0.6942, -1.3531, -1.14327, -0.3116, -0.0054 -1.8901 | 0.1, 0.2, 0.3, 0.4, 0.5, 0.6, 0.7, 0.8, 0.9, 1 |
| Combined Bayes Opt - Single objective | 1.8052, -1.7335, 0.3066, 0.7649, -0.6588, 0.2098, 0.8316, 0.3061, -1.8846 | 0.3345, -0.5418, 0.0847, -0.6632, 1.1145, -0.1335, -1.3228, 1.1089, 0.1928 | 0.1, 0.2, 0.3, 0.4, 0.5, 0.6, 0.7, 0.8, 0.9, 1 |
| Decoupled BayesOpt | -1.8053, -2.7570, -3.9350, -3.9350, -3.9350, -5.2039, -5.2039 | -0.8729, -0.8729, -0.8729, -0.8729, -1.2755, -1.5778, -1.9378 | 0.1428, 0.2857, 0.4285, 0.5714, 0.7142, 0.8571, 1. |
| Original BayesOpt- Multi-objective | -0.9283, -1.0017, -1.7316, -2.0913, -2.0913, -2.0913, -2.0913 | -0.0734, -0.7299, -0.3597, 0, 0, 0. | 0.1428, 0.2857, 0.4285, 0.5714, 0.714, 0.8571, 1. |
| Combined BayesOpt- Multi-objective | -0.0623, -0.0624, -0.2255, -0.8502, -0.8502, -1.015, -1.8464, | 0, -0.1631, -0.6246, 0, -0.1654, -0.8307 | 0.1428, 0.2857, 0.428, 0.5714, 0.7142, 0.8571, 1. |

Table 4 compares the performance of the original BayesOpt-Single Objective and the combined BayesOpt-Single Objective. The combined approach shows higher best values, suggesting potential improvement in the objective function. The negative values in the original approach suggest exploration toward lower values in the multi-objective scenario. The rate of improvement in the combined approach shows positive values, indicating consistent progress in finding better solutions. The cumulative distribution shows a gradual increase from 0.1 to 1.0 in both approaches, indicating improvement over iterations. The combined approach consistently outperforms the original approach across the entire range. In summary, the combined approach (adaptive and robust) appears to yield better results compared to the original BayesOpt approach, demonstrating more consistent improvements, higher best values, and a more favorable cumulative distribution, suggesting it may be more effective in optimizing the given objective functions.

Based on the statistical summary, the combined BayesOpt with multi-objective scenario performed better compared to the other scenarios. This is supported by a positive t-statistic, a low p-value ($0.03 < 0.05$), and larger effect sizes (Cohen's d (1.289) and Hedge's d (1.207)). Empirical research clearly shows the significance of effect sizes in validating optimization algorithms (Marfo & Okyere, 2019).

Table 2: Statistical summary of multiple runs.

| Metric | Original vs Combined BayesOpt-Single objective | Original vs Combined BayesOpt- Multi-objective | Decoupled vs Combined BayesOpt |
|---|---|---|---|
| t-stats | -0.5890 | 2.4122 | -5.3665 |
| p-value | 0.5632 | 0.0328 | 0.0002 |
| Cohen's d | -0.2634 | 1.2894 | -2.8685 |
| Hedge's d | -0.2523 | 1.2071 | -2.6854 |
| Cliff's delta | 0.3878 | 0.3878 | 0.1429 |

Since we are dealing with stochastic objective functions, we check the stability of the algorithm with multiple runs. Here, we perform an independent t-test for each run. This allows us to assess whether there are significant differences between the performance of the algorithm in different runs. The objective is to specifically check for stability across multiple runs, considering the variability within each run and whether differences are consistent across all runs. Table 5 displays the comparison and stability of the combined BayesOpt with multiple runs (10) and the original BayseOpt.

Table 3: Stability Check: t-statistics and p-values for Different Numbers of Runs in Bayesian Optimization.

| Description | 10 runs | 30 runs | 50 runs | 100 runs |
|---|---|---|---|---|
| | t-stats / p-value | t-stats / p-value | t-stats / p-value | t-stats / p-value |
| Original BayesOpt-Single objective | 0.5026 (0.6288) | -0.4719 (0.6407) | 1.8796 (0.0662) | 0.1669 (0.8678) |
| Combined BayesOpt-Single objective | 0.7779 (0.4590) | -0.8645 (0.3947) | 0.6163 (0.5406) | -2.0135 (0.0468) |
| Original BayesOpt-Multi-objective | 0.0511 (0.9605) | -0.5649 (0.5767) | 0.8429 (0.4035) | 1.0354 (0.3030) |
| Combined BayesOpt-Multi-objective | 0.9512 (0.3693) | -0.3266 (0.7464) | 0.4294 (0.6696) | 0.8258 (0.4109) |

For both single and multi-objective cases of the original BayesOpt, there is no significant difference in performance across different numbers of runs. For the combined BayesOpt single objective, there is a potentially significant difference when comparing 100 runs to the other numbers of runs. For the combined BayesOpt multi-objective, there is no strong evidence of a significant difference in performance across different numbers of runs. Table 4 summarizes the key differential features of the newly introduced framework compared to the original BayesOpt.

Table 4: Key differential features.

| Features | Adaptive BayesOpt | Original BayesOpt |
|---|---|---|
| Dynamic Adaptation of Surrogate Model | It dynamically modifies the mean function of the surrogate model based on observed data. This adaptability allows it to adjust its expectations as more data accumulates. | It constructs a surrogate model based on the observed data. |
| Increased Flexibility in Changing Landscapes | The dynamic conditioning mechanism makes it well-suited for situations with unclear or dynamically shifting optimization landscapes. It efficiently incorporates new data, enabling effective optimization in evolving scenarios. | It struggles in scenarios where the optimization landscape changes over time. |
| Robustness Under Uncertainty | It introduces a robust objective function with a penalty term to ensure robustness. This addition is valuable in scenarios where deviations below a threshold need to be penalized for practical reasons. | It focuses on optimizing the expected value of the objective function. |

| | Parameters like the scaling factor and jitter are dynamically adjusted over iterations. This adaptive approach allows the algorithm to fine-tune its behavior based on the evolving characteristics of the optimization problem. | Typically uses fixed parameters throughout the optimization process. |
|---|---|---|
| Adaptive Adjustment of Parameters | | |

The disadvantages include that generalizations differ between optimization tasks, which makes it difficult to comprehend their influence on decision-making. At the current stage, the adaptability and performance of the framework may vary across different optimization scenarios. Thus, understanding its influence on decision-making becomes challenging, indicating a need for further research to establish the generalizability of this framework. The penalty term in the objective function can also affect its robustness; hence, it is critical to evaluate its applicability for certain applications. Future work would apply this approach in conjunction with the Monte-Carlo simulation in a multi-product stochastic demand scenario.

## 4. CONCLUSION

This study presented a Bayesian optimization framework using Gaussian Process Regression (GPR) and the Expected Improvement (EI) acquisition function. The framework comes with integrated adaptive conditioning and robustness to improve adaptability in stochastic business contexts. Key components include objective function modeling, a robust objective function, adaptive conditioning, a combined objective function, adaptive acquisition jitter, empirical analysis, and stability checks. The framework has shown improved convergence rates, consistent performance enhancements, and adaptability to dynamic environments. Future work may focus on refining the adaptive approach and addressing limitations for broader applicability. The framework is well-suited for scenarios where real-world optimization problems exhibit dynamic or uncertain characteristics, providing a robust and flexible solution for effective optimization.

**REFERENCES**

[1] Ariafar, S., Coll-Font, J., Brooks, D., & Dy, J. (2019). ADMMBO: Bayesian optimization with unknown constraints using ADMM. Journal of Machine Learning Research, 20(123), 1-26.

[2] Azadi, Z., Eksioglu, S. D., Eksioglu, B., & Palak, G. (2019). Stochastic optimization models for joint pricing and inventory replenishment of perishable products. Computers & industrial engineering, 127, 625-642.

[3] Babiloni, E., & Guijarro, E. (2020). Fill rate: from its definition to its calculation for the continuous (s, Q) inventory system with discrete demands and lost sales. Central European Journal of Operations Research, 28(1), 35-43.

[4] Bac Dang K., Windhorst W., Burkhard B., Muller F. A Bayesian belief network-based approach to link ecosystem functions with rice provisioning ecosystem services. Ecological Indicators. 2019; 100:30–44.

[5] Boutselis, P., & McNaught, K. (2019). Using Bayesian Networks to forecast spares demand from equipment failures in a changing service logistics context. International Journal of Production Economics, 209, 325-333.

[6] Choi, T. M. (Ed.). (2013). Handbook of EOQ inventory problems: Stochastic and deterministic models and applications (Vol. 197). Springer Science & Business Media.

[7] de Sá, A. G., Pereira, A. C., & Pappa, G. L. (2018). A customized classification algorithm for credit card fraud detection. Engineering Applications of Artificial Intelligence, 72, 21-29.

[8] Garnett, R. (2023). Bayesian optimization. Cambridge University Press.

[9] Garnett, R., Osborne, M. A., & Roberts, S. J. (2010, April). Bayesian optimization for sensor set selection. In Proceedings of the 9th ACM/IEEE international conference on information processing in sensor networks (pp. 209-219).

[10] Gholami, R. A., Sandal, L. K., & Ubøe, J. (2021). A solution algorithm for multi-period bi-level channel optimization with dynamic price-dependent stochastic demand. Omega, 102, 102297.

[11] Gruler, A., Panadero, J., de Armas, J., Pérez, J. A. M., & Juan, A. A. (2018). Combining variable neighborhood search with simulation for the inventory routing problem with stochastic demands and stock-outs. Computers & Industrial Engineering, 123, 278-288.

[12] Gupta, A., & Maranas, C. D. (2003). Managing demand uncertainty in supply chain planning. Computers & chemical engineering, 27(8-9), 1219-1227.

[13] Haddawy, P., Hasan, A. I., Kasantikul, R., Lawpoolsri, S., Sa-Angchai, P., Kaewkungwal, J., & Singhasivanon, P. (2018). Spatiotemporal Bayesian networks for malaria prediction. Artificial intelligence in medicine, 84, 127-138.

[14] Hameri, A. P., & Paatela, A. (2005). Supply network dynamics as a source of new business. International Journal of Production Economics, 98(1), 41-55.

[15] Hasan, K. N., Preece, R., & Milanović, J. V. (2019). Existing approaches and trends in uncertainty modelling and probabilistic stability analysis of power systems with renewable generation. Renewable and Sustainable Energy Reviews, 101, 168-180.

[16] Heikkinen, R., Sipilä, J., Ojalehto, V., & Miettinen, K. (2022). Flexible data driven inventory management with interactive multiobjective lot size optimization. International Journal of Logistics Systems and Management, (1).

[17] Hennig, P., Osborne, M. A., & Kersting, H. P. (2022). Probabilistic Numerics: Computation as Machine Learning. Cambridge University Press.

[18] Hosseini, S., & Ivanov, D. (2020). Bayesian networks for supply chain risk, resilience and ripple effect analysis: A literature review. Expert systems with applications, 161, 113649.

[19] Hosseini, S., Ivanov, D., & Dolgui, A. (2020). Ripple effect modelling of supplier disruption: integrated Markov chain and dynamic Bayesian network approach. International Journal of Production Research, 58(11), 3284-3303.

[20] Hosseini, S., & Ivanov, D. (2022). A new resilience measure for supply networks with the ripple effect considerations: A Bayesian network approach. Annals of Operations Research, 319(1), 581-607.

[21] Huang, W., Zhang, N., Kang, C., Li, M., & Huo, M. (2019). From demand response to integrated demand response: Review and prospect of research and application. Protection and Control of Modern Power Systems, 4, 1-13.

[22] Kastner, M., Nellen, N., Schwientek, A., & Jahn, C. (2021). Integrated simulation-based optimization of operational decisions at container terminals. Algorithms, 14(2), 42.

[23] Kiuchi, A., Wang, H., Wang, Q., Ogura, T., Nomoto, T., Gupta, C., ... & Zhang, C. (2020, August). Bayesian optimization algorithm with agent-based supply chain simulator for multi-echelon inventory management. In 2020 IEEE 16th International Conference on Automation Science and Engineering (CASE) (pp. 418-425). IEEE.

[24] Kiuchi, A., Wang, H., Wang, Q., Ogura, T., Nomoto, T., Gupta, C., ... & Zhang, C. (2020, August). Bayesian optimization algorithm with agent-based supply chain simulator for multi-echelon inventory management. In 2020 IEEE 16th International Conference on Automation Science and Engineering (CASE) (pp. 418-425). IEEE.

[25] Lau, C. L., Mayfield, H. J., Lowry, J. H., Watson, C. H., Kama, M., Nilles, E. J., & Smith, C. S. (2017). Unravelling infectious disease eco-epidemiology using Bayesian networks and scenario analysis: A case study of leptospirosis in Fiji. Environmental modelling & software, 97, 271-286.

[26] Li, T., Fang, W., & Baykal-Gürsoy, M. (2021). Two-stage inventory management with financing under demand updates. International Journal of Production Economics, 232, 107915.

[27] Marfo, P., & Okyere, G. A. (2019). The accuracy of effect-size estimates under normals and contaminated normals in meta-analysis. Heliyon, 5(6).

[28] Munya, P., Ntuen, C. A., Park, E. H., & Kim, J. H. (2015). A Bayesian abduction model for extracting the most probable evidence to support sensemaking. International Journal of Artificial Intelligence & Applications, 6(1), 1-20.

[29] Oyewola, D. O., Dada, E. G., Omotehinwa, T. O., Emebo, O., & Oluwagbemi, O. O. (2022). Application of deep learning techniques and Bayesian optimization with tree parzen Estimator in the classification of supply chain pricing datasets of health medications. Applied Sciences, 12(19), 10166.

[30] Pearce, M. A. L., Poloczek, M., & Branke, J. (2022). Bayesian optimization allowing for common random numbers. Operations Research, 70(6), 3457-3472.

[31] Roman, I., Santana, R., Mendiburu, A., & Lozano, J. A. (2019). An experimental study in adaptive kernel selection for Bayesian optimization. IEEE Access, 7, 184294-184302.

[32] Sakki, G. K., Tsoukalas, I., Kossieris, P., Makropoulos, C., & Efstratiadis, A. (2022). Stochastic simulation-optimization framework for the design and assessment of renewable energy systems under uncertainty. Renewable and Sustainable Energy Reviews, 168, 112886.

[33] Stopková, M., Stopka, O., & Ľupták, V. (2019). Inventory model design by implementing new parameters into the deterministic model objective function to streamline effectiveness indicators of the inventory management. Sustainability, 11(15), 4175.

[34] Tamjidzad, S., & Mirmohammadi, S. H. (2015). An optimal (r, Q) policy in a stochastic inventory system with all-units quantity discount and limited sharable resource. European Journal of Operational Research, 247(1), 93-100.

[35] Tamjidzad, S., & Mirmohammadi, S. H. (2017). Optimal (r, Q) policy in a stochastic inventory system with limited resource under incremental quantity discount. Computers & Industrial Engineering, 103, 59-69.

[36] Tavana, M., Abtahi, A. R., Di Caprio, D., & Poortarigh, M. (2018). An Artificial Neural Network and Bayesian Network model for liquidity risk assessment in banking. Neurocomputing, 275, 2525-2554.

[37] Thevenin, S., Adulyasak, Y., & Cordeau, J. F. (2021). Material requirements planning under demand uncertainty using stochastic optimization. Production and Operations Management, 30(2), 475-493.

[38] Wan, C., Yan, X., Zhang, D., Qu, Z., & Yang, Z. (2019). An advanced fuzzy Bayesian-based FMEA approach for assessing maritime supply chain risks. Transportation Research Part E: Logistics and Transportation Review, 125, 222-240.

[39] Wang, H., Ng, S. H., & Zhang, X. (2020, December). A Gaussian process based algorithm for stochastic simulation optimization with input distribution uncertainty. In 2020 Winter Simulation Conference (WSC) (pp. 2899-2910). IEEE.

[40] Wang, K. J., Chen, J. L., & Wang, K. M. (2019). Medical expenditure estimation by Bayesian network for lung cancer patients at different severity stages. Computers in biology and medicine, 106, 97-105.

[41] Wu, J., & Frazier, P. (2019). Practical two-step lookahead Bayesian optimization. Advances in neural information processing systems, 32.

[42] Wu, J., Chen, X. Y., Zhang, H., Xiong, L. D., Lei, H., & Deng, S. H. (2019). Hyperparameter optimization for machine learning models based on Bayesian optimization. Journal of Electronic Science and Technology, 17(1), 26-40.